\documentclass[12pt]{paper}
 \usepackage{centernot}
 \usepackage{epsfig}
\usepackage{graphics}
 \usepackage{color}
 \newcommand{\tug}{\tilde u_g}
 \newcommand{\twg}{\tilde w_g}
 \newcommand{\hug}{\hat u_g}
  \newcommand{\hwg}{\hat w_g}
 \newcommand{\Ho}{H^1_0([0,T];L^2(\Gamma))}
 \newcommand{\No}{{\mathcal{N}_1}}
 
 \newcommand{\Nod}{{\mathcal{N}_2}}
  
\newcommand{\Mo}{\mathbb{M}}

\usepackage{color}
\usepackage{epsfig,cmmib57}
\usepackage{amssymb,amsmath,exscale}
 \numberwithin{equation}{section}

\newcommand{\ZEP}{\epsilon}

\newcommand{\ZSUno}{\sum _{n=1}^{+\ZIN}}
\newcommand{\ZOMq}{\Omega}

\newcommand{\ZA}{{\mathcal A}}

\newcommand{\zg}{\gamma}

\newcommand{\intT}{\int_0^T}
\newcommand{\intt}{\int_0^t}
\newcommand{\ints}{\int_0^s}
\newcommand{\intr}{\int_0^r}

\newcommand{\ZCD}{(\cdot)}
\newtheorem{Theorem}{Theorem}  

\newtheorem{Lemma}[Theorem]{Lemma}

\newtheorem{Remark}[Theorem]{Remark}
\newtheorem{Definition}[Theorem]{Definition}

\newcommand{\zdiaform}{\mbox{~~\zdia}}

\newcommand{\zaa}{\alpha}

\newcommand{\ZDE}{\delta}
\newcommand{\zt}{\tau}
\newcommand{\zdia}{~~\rule{1mm}{2mm}\par\medskip}

\newcommand{\ZLA}{\label}
\newcommand{\ZIN}{\infty}
\newcommand{\zProof}{{\noindent\bf\underbar{Proof}.}\ }

\newcommand{\zzr}{{\rm I\hskip-2.1pt R}}
\newcommand{\ZBI}{\bibitem}
\newcommand{\ZD}{\;\mbox{\rm d}}
\newcommand{\zl}{\lambda}

\author{
L. Pandolfi\thanks{Dipartimento di Scienze Matematiche ``Giuseppe Luigi Lagrange'', Politecnico di Torino, Corso Duca degli Abruzzi 24, 10129 Torino, Italy (luciano.pandolfi@polito.it)}
}
\title{Controllability and lack of controllability with smooth controls in viscoelasticity via moment methods\thanks{
This papers fits into the research program of the GNAMPA-INDAM and has been written in the framework of the   ``Groupement de Recherche en Contr\^ole des EDP entre la France et l'Italie (CONEDP-CNRS)''.}}
\begin{document}
 
 \maketitle 
 
\noindent {\bf\underline{Abstract}:} In this paper we study controllability of a linear  equation with persistent memory when the control belongs to $ H^k_0(0,T;L^2(\ZOMq)) $. 
In the case the memory is zero, our equation is reduced to the wave equation and a result
due to Everdoza and Zuazua  informally states  that smoother targets can be reached by using smoother controls.  In this paper we prove that this result can be partially extended to systems with memory, but that the memory is an obstruction to a complete extensions.

\medskip


\smallskip


\section{Introduction}
Stimulated by the applications to the quadratic regulator problem, controllability for distributed parameter systems is usually studied with square integrable controls.  Such general controls 
   are hardly realizable in practice and only smooth or piecewise smooth controls,  like bang-bang controls, can be implemented. 
Moreover, when a control is implemented numerically, via discretization, convergence estimates depend on the smoothness of the control  (see~\cite{EverdZUAZlibro}).    
     This fact revived interest on controllability with smooth controls, and the crucial results are in~\cite{EverdZUAZArch}     
(see also~\cite[Theorem~2.1]{Russel} where reachability of the wave equation under controls of class $H^s$  acting on the entire boundary  is studied).
A natural guess (which is true for the wave equation but which we in part disprove for systems with memory) is  that if  the control is smooth then the reachable targets are ``smooth'' too, and the problem is to identify the targets which can be reached by using controls in certain smoothness classes.

In this paper we are going to examine 
  the following equation, which is encountered in viscoelasticity and in diffusion processes when the material  has a complex molecular structure:
\begin{equation}
\ZLA{eq:modello}
w''=\left (\Delta w+bw\right )+\intt K(t-s)w(s)\ZD s \,.
\end{equation}
Here $w=w(x,t)$, the apex denotes time derivative, $w''(x,t)=w_{tt}(x,t)$,  $x\in\ZOMq\subseteq \zzr^d$ is a bounded region with $C^2$ boundary, $K(t)$ is a real continuous function  and $\Delta=\Delta_x$ is the laplacian in the variable $x$.

Dependence on the time and expecially space variable is not explicitly indicated unless needed for clarity so that we shall write $w=w(t)=w(x,t)$ according to convenience.

We associate the following initial/boundary conditions to system~(\ref{eq:modello}):
\begin{equation}
\ZLA{eq:iniBOUNDcond}
\begin{array}
{l}
w(0)=w_0\,,\quad w'(0)=w_1\,,\\
w(x,t)= \left\{\begin{array}{cc} f(x,t) &x\in \Gamma\\
0 & x\in\partial\ZOMq\setminus\Gamma 
\end{array}\right.
\end{array}
\end{equation}
($\Gamma $ is a relatively open subset of $\partial\ZOMq$).

The function $f$ is a control, which is used to steer the pair $\left (w(t),w'(t)\right )$ to hit a prescribed target $\left  (\xi,\eta\right )$ at a certain time $T$.  

The spaces of the initial data and final targets and of the control $f$ will be specified below.

There is no assumption on the sign of $b$ whose presence is explained in Remark~\ref{RemaMacC}. Furtermore we note:

\begin{itemize}
\item It is known (and recalled in Sect.~\ref{subsAssumpDISCUSS})
 that when $f\in L^2(0,T;L^2(\Gamma))$ and $( w_0,w_1)\in L^2(\ZOMq)\times H^{-1}(\ZOMq)$ then problem~(\ref{eq:modello})-(\ref{eq:iniBOUNDcond}) admits a unique mild solution $(w(t),w'(t))\in C([0,T];L^2(\ZOMq)\times H^{-1}(\ZOMq))$ for every $T>0$;
 \item when $K=0$ (i.e. when we consider the wave equation) the solution of~(\ref{eq:modello}) is denoted $u$; 
\item when we want to stress the dependence on $f$ of the solution of~(\ref{eq:modello}) we use the notation $w_f$ (the notation $u_f$ when $K=0$).\zdia
  \end{itemize}

In order to describe the result proved in~\cite{EverdZUAZArch} it is convenient to introduce the following operators $A$ and $\ZA$ in $L^2(\ZOMq)$:
\begin{equation}\ZLA{eq:defiOpeA}
{\rm dom}\, A=H^2(\ZOMq)\cap H^1_0(\ZOMq)\,,\qquad A\phi=\Delta\phi+b\phi\,,\quad \ZA=(-A)^{1/2}
\end{equation}
(note that $A$ is a positive operator if $b\geq 0$ and in this case $\ZA$ is defined in a 
standard way; if $b<0$ the definition of $\ZA$ is discussed in Sect.~\ref{subsAssumpDISCUSS}).

It turns out that ${\rm dom}\,\ZA=H^1_0(\ZOMq)$.

\begin{Definition}
Let $T>0$ and let $\mathcal{F}$ be a closed subspace of $ L^2(0,T;L^2(\Gamma))$.
 We say that Eq.~(\ref{eq:modello}) is ${\rm dom}\,\ZA^{k+1}\times {\rm dom}\,\ZA^k$-controllable at time $T$ with controls $f\in \mathcal{F}$ when the following properties hold:
\begin{enumerate}
 \item
If $(w_0,w_1)\in {\rm dom}\,\ZA^{k+1}\times {\rm dom}\,\ZA^k    $  and $f\in \mathcal{F}$ then $(w_f(T),w_f'(T))\in  {\rm dom}\,\ZA^{k+1}\times {\rm dom}\,\ZA^k $;
 \item
 for every $w_0$, $\xi_0$ in $ {\rm dom}\,\ZA^{k+1}$ and every $w_1$, $\eta$ in ${\rm dom}\,\ZA^k$ there exists $f\in \mathcal{F}$ such that
 $(w_f(T),w_f'(T))=(\xi,\eta)$. 

\end{enumerate}
\end{Definition}
%It is easily seen that when studying controllability we can assume $w_0=0$ and $w_1=0$.

The following result is proved   in~\cite{PandSharp} (see also~(\cite{PandLibro,PandParma,PandESAIM}).
\begin{Theorem}\ZLA{Theo:controLdue}
  There exists a time $T_0$ and a relatively open set $\Gamma\subseteq\partial\ZOMq$ 
which have the following properties: Let $T>T_0$. For every $w_0$ and $\xi$ in $L^2(\ZOMq) $ and for every $w_1$ and $\eta$ in $H^{-1}(\ZOMq)$ 
  there exists $f\in L^2(0,T;L^2(\partial \ZOMq))$ such that $\left (w_{ f}(T),w'_{  f}(T)\right )=\left (\xi,\eta\right )$.
  The set $ \Gamma $ and the number $ T_0 $ do not depend on the continuous memory kernel $ K(t) $.
\end{Theorem}
Note that Theorem~\ref{Theo:controLdue} holds in particular for the wave equation (i.e. when $K=0$) and the proofs in the references above do  depend on the known controllability result of the wave equation.
 
The result in~\cite{EverdZUAZArch} can be adapted to the case of the wave equation (without memory) as described in~\cite[Sect.~5.2]{EverdZUAZArch}  (see item~\ref{item1inRemaMacC} in Remark~\ref{RemaMacC} to understand the exponents):

\begin{Theorem}\ZLA{TeoDIevErdZUazPERiperb}
Let $T$, $T_0$ and $\Gamma$ be as in Theorem~\ref{Theo:controLdue}. System~(\ref{eq:modello}) with $K=0$ is
${\rm dom}\,\ZA^{k }\times {\rm dom}\,\ZA^{k-1}$-controllable at time $T$ with controls  $f\in H^k_0(0,T;L^2(\Gamma))$.
\end{Theorem}
In the light of Theorem~\ref{Theo:controLdue} (which extends the well known controllability result of the wave equation) it is natural to guess that Theorem~\ref{TeoDIevErdZUazPERiperb} can be extended too. Instead we have the following result:

\begin{Theorem}\ZLA{TeoINh10} 
Let $T_0$, $ T $ and $\Gamma$ have the properties in Theorem~\ref{Theo:controLdue}.  Then we have:
\begin{enumerate}
\item
\ZLA{item1delTeoINh10}  
 
System~(\ref{eq:modello}) is
${\rm dom}\,\ZA \times L^2(\ZOMq)$-controllable at time $T$ with controls  $f\in H^1_0(0,T;L^2(\Gamma))$ (note that $L^2(\ZOMq)={\rm dom}\,\ZA^0$).
 \item\ZLA{item2delTeoINh10} 
System~(\ref{eq:modello})  
${\rm dom}\,\ZA^{2}\times {\rm dom}\,\ZA $-controllable at time $T$ with controls  $f\in H^2_0(0,T;L^2(\Gamma))$.

\item
\ZLA{item3delTeoINh10}    
Let $k\geq 3$. For every $T>0$ there exist controls  $f\in H^k_0(0,T;L^2(\Gamma))$ such that $(w(T),w'(T))\notin
{\rm dom}\,\ZA^{k}\times {\rm dom}\,\ZA^{k-1}  $.
  
  \end{enumerate}
\end{Theorem}

\begin{Remark}\ZLA{RemaMacC}
{\rm
\begin{enumerate}
\item\ZLA{item1inRemaMacC} The operator $A$ in~\cite[Sect.~5.2]{EverdZUAZArch} is defined as the laplacian with domain $H^1_0(\ZOMq)$ (and image $H^{-1}(\ZOMq)$ while we used ${\rm dom}\, A=H^2(\ZOMq)\cap H^1_0(\ZOMq)$.  
 \item for the sake of brevity, the properties in Theorems~\ref{Theo:controLdue}   will be called   ``controllability in $L^2(\ZOMq)\times H^{-1}(\ZOMq)$''.
 \item In the case $\ZOMq=(a,b)$, item~\ref{item1delTeoINh10} of Theorem~\ref{TeoINh10}
 has been proved in~\cite{PandTRIULZI}.
\item it is clear that when studying controllability we can reduce ourselves to study the system with zero initial conditions, $w_0=0$, $w_1=0$.
 
\item
The usual form of the system with persistent memory which is encountered in viscoelasticity is
\[
w''=\Delta w+\intt M(t-s)\Delta w(s)\ZD s\,.
\]
We formally solve this equation as a Volterra integral equation in the ``unknown'' $\Delta w$. Two     integrations  by parts in time (followed by an exponential transformation) lead  to Eq.~(\ref{eq:modello}),  with $b\neq 0$ (if the initial conditions are different from zero also an affine term, which contains the initial conditions appear, but when studying controllability we 
can assume $w_0=0$, $w_1=0$). For this reason we kept the addendum $bw$ in Eq.~(\ref{eq:modello}). This transformation is known as \emph{MacCamy trick} and it is detailed in~\cite{PandLibro}.\zdia 
 
\end{enumerate}
}
\end{Remark}

% It is convenient to define the reachable spaces:
%\begin{align*}
%%R_E(T)&=\left \{ \left (u_f(T),u_f'(T)\right )\,,\ f\in L^2(0,T;U)\right\}\,,\\
% R_V(T)&=\left \{ \left (w_f(T),w_f'(T)\right )\,,\ f\in L^2(0,T;U)\right \}\,,\\
%%   R_{1, E}(T)&=\left \{ \left (u_f(T),u_f'(T)\right )\,,\ f\in H^1_0(0,T;U)\right \}\,,\\
%  R_{k,V}(T)&=\left \{ \left (w_f(T),w_f'(T)\right )\,,\ f\in H^k_0(0,T;U)\right \}\,.
%\end{align*}
% When $K=0$, i.e. in the case of the wave equation, these spaces are denoted $R_{k,E}(T)$.
%
%\CM{Si usa davvero questa notazione?}

\subsection{\ZLA{subsAssumpDISCUSS}Preliminaries}
 
The operator $A$     in $L^2(\ZOMq)$ was already defined: $A\phi=\Delta\phi+b\phi$, 
${\rm dom}\, A=H^2(\ZOMq)\cap H^1_0(\ZOMq)$ (we recall that $\ZOMq$ is a region with $C^2$ boundary).

The operator $A$ is selfadjoint, possibly non positive since 
  there is no assumption on the sign of $b$. Its resolvent is compact so that  the Hilbert space $L^2(\ZOMq)$ has  an orthonormal basis $\{\phi_n\}$ of eigenvectors   of the operator $A$. We denote $-\zl_n^2$ the eigenvalue of $\phi_n$ since $ \zl_n^2>0$   for large $n$ (it might be $ \zl_n^2\leq 0$ if $n$ is small).    The eigenvalues are repeated according to their multiplicity (which is finite).
 
 We shall use the following known asymptotic estimate for the eigenvalues (see~\cite{Agmon}):
\begin{equation}\ZLA{eq:stimAUTOV}
\mbox{if $\ZOMq\subseteq\zzr^d$ then $\zl_n^2\sim n^{2/d}$}\,.
\end{equation}  
In particular, $\zl_n^2$
   is positive for large $n$.
   
Let $\gamma\in(0,1)$. If $n$ is large and $\zl_n^2\geq 0$ then $\zl_n^{2\zg}$ is the nonegative determination; otherwise we fix one of the   determinations.  When $\zg=1/2$,  for example the one with nonnegative imaginary part. \emph{So, $\zl_n $ denotes the chosen determination of the square root of  $\zl_n^2$, $ \zl_n>0 $ is $ n $ is large.}

By definition,

\[
\begin{array}{l}
\displaystyle \xi=\ZSUno \xi_n\phi_n(x)\in {\rm dom }\,(-A)^\zg  \ \iff
\{\zl_n^{2\zg}\xi_n\}\in l^2 \\[2mm]
\displaystyle  {\rm and}\quad 
 (-A)^\zg\xi= \ZSUno \zl_n^{2\zg}\xi_n\phi_n \quad   \mbox{so that in particular}\\
\ZA\xi =i\left (-A\right )^{1/2}\xi=i\left (\ZSUno \zl_n\xi_n\phi_n(x)\right ) 
  \,.
 \end{array}
\]
 Furthermore we define (we recall the that $\zl_n$ is real when $n$ is large)
 \begin{align*}
&R_+(t)\left ( \ZSUno \xi_n\phi_n(x)\right )= \ZSUno \left (\cos\zl_n t\right ) \xi_n\phi_n(x)\,, 
 \\
 &
R_-(t)\left ( \ZSUno \xi_n\phi_n(x)\right )= i\left (\ZSUno \left (\sin \zl_n t\right ) \xi_n\phi_n(x)\right ) \,.
 \end{align*}

 Finally, we introduce the operator $D$: $L^2(\Gamma)\mapsto L^2(\ZOMq)$:
\[
u=Df\ \iff \ \left\{\begin{array}{l}
\Delta u+bu=0\ {\rm in}\ \ZOMq\,,\\ u=  f\ {\rm on}\, \Gamma\,,\quad u=0\ {\rm on}\  \partial\ZOMq\setminus\Gamma\,.
\end{array}\right.
\]
It is known that  ${\rm im}\, D\subseteq H^{1/2}(\ZOMq)\subseteq {\rm dom}\, (-A)^{1/4-\ZEP} $ for every $\ZEP>0$. 
 
It is known (see~\cite{LasTRIEquonde}) that the mild solution of  the wave equation
\[
u''=\Delta u+F\,, 
\] 
with initial and boundary conditions~(\ref{eq:iniBOUNDcond}) is 
 \begin{align}\nonumber
u(t)=&R_+(t)w_0+\ZA^{-1}R_-(t)w_1-\ZA\intt R_-(t-s) D f(s)\ZD s\\
\ZLA{eqDIu}
&+\ZA^{-1}\intt R_-(t-s)F(s)\ZD s\,.  
  \end{align}

By definition, the mild  solution of problem~(\ref{eq:modello})-(\ref{eq:iniBOUNDcond})  
  is the solution of the following Volterra integral equation
  \begin{equation}
  \ZLA{eq:soluVisco}
w (t)=u (t)+\ZA^{-1}\intt\left [\int_0^{t-s} K(r)R_-(t-s-r) w(s)\ZD r\right ] \ZD s 
  \end{equation}
 where $u(t)$ is given by~(\ref{eqDIu}) with $F=0$.

    We note that $w'(t)$ is given by
\begin{equation}
  \ZLA{eq:volteDERIVATA}
  w'(t)=u'(t) +\intt \left [\int_0^{t-s} K(s)R_+(t-s-r) w(s)\ZD r\right ]\ZD s
  \end{equation}
where
\begin{equation}
\ZLA{eq:diuPRIMO}
u'(t)=\ZA R_-(t)w_0+R_+(t) w_1-A\intt R_+(t-s) D f(s)\ZD s\,.
\end{equation}
The following result is known (see~\cite[Ch.~2]{PandLibro}):
\begin{Theorem}\ZLA{teo:propriESoluOnde}
If $(\xi,\eta,f)\in L^2(\ZOMq)\times H^{-1}  (\ZOMq)\times  L^2\left (0,T;L^2(\partial\ZOMq)\right )$ 
then $\left (w_f(t),w_f'(t)\right )\in C\left ([0,T];L^2(\ZOMq)\right )\times C\left ([0,T];H^{-1}(\ZOMq)\right )$
 for every $T>0$ and the linear transformation 
$(\xi,\eta,f)\mapsto (w,w')$ is continuous in the indicated spaces.
\end{Theorem}

Finally, let $w_0=0$, $w_1=0$. We introduce the following operators  which are continuous from $ L^2(0,T;L^2(\Gamma)) $ to $ L^2(\ZOMq)\times H^{-1}(\ZOMq) $:
 
\begin{align*}
&f\mapsto \Lambda_E(T) f=\left ( \begin{array}{cc}
\Lambda _{ E, 1}(T) f&\Lambda _{E, 2}(T)f
\end{array}\right)= \left ( \begin{array}{cc}
u_ f(T)&u'_f(T)
\end{array}\right)\\
&f\mapsto \Lambda_V(T) f=\left ( \begin{array}{cc}
\Lambda _{ V,1} (T)f&\Lambda _{ V,2}(T)f
\end{array}\right)= \left ( \begin{array}{cc}
w_ f(T)&w'_f(T)
\end{array}\right)\,.
\end{align*}
\emph{Our assumption is that the set $\Gamma$ and the time $T$ have been chosen so that both these operators are surjective.}
%Note that 
%\[
%%R_E(T)={\rm im}\, \Lambda_E(T)\,,\qquad 
% R_V(T)={\rm im}\, \Lambda_V(T)   \,,\qquad R_{1,V}(T)=\Lambda_V(T)_{|_{H^1_0(0,T;L^2(\ZOMq))}}\,.
%\]
%while
%\[
%R_{1,V}(T) \quad \mbox{is the image of ${\rm im}\, \Lambda_V(T)  $ restricted to $ H^1_0(0,T;L^2(\Gamma)) $.}
%\]
%%\emph{Our assumption is that the set $\Gamma$ and the time $T$ have been chosen so that both these operators are surjective.}
%\CM{Credo che queste notazioni non si usino}
\section{Controllability with square integrable controls and moment problem}
From now on we study controllability   and so we assume $w_0=0$, $w_1=0$.

We expand the solutions of Eq.~(\ref{eq:modello}) in series of $\phi_n$,
\begin{equation}\ZLA{eq:expaWn}
w(x,t)=\ZSUno \phi_n(x) w_n(t)\,,\qquad w_t(x,t)=\ZSUno \phi_n(x) w_n'(t)\,.
\end{equation}
It is easily seen that $w_n(t)$ solves 
\[%  
w_n''(t)=-\zl_n^2 w_n(t)+\intt K(t-s)w_n(s)\ZD s-\int_\Gamma \left (\zg_1\phi_n\right ) f(x,t)\ZD \Gamma
\]
where $\zg_1$ is the exterior normal derivative and $\ZD\Gamma$ is the surface measure.

The initial conditions are zero since $w(0)=0$, $w'(0)=0$.
  In order to represent the solution of the previous equation, we introduce $\zeta_n(t)$, the solution of      
\begin{equation}
\ZLA{eq:diZETAn}
\zeta_n''(t)=-\zl_n^2\zeta_n(t)+\intt K(t-s) \zeta_n(s)\ZD s\,,\qquad \left\{\begin{array}{l}
\zeta_n(0)=0\,,\\
\zeta_n'(0)=1\,.
\end{array}\right.
\end{equation}
 Then we have     
\begin{equation}\ZLA{RApprewN}
\left\{
\begin{array}{l}
\displaystyle w_n(t)=\int_\Gamma\intt \left [\zeta_n(t-s)\zg_1\phi_n(x)\right ] f(x,s)\ZD s\,,\\[3mm]
\displaystyle w_n'(t)=\int_\Gamma\intt \left [\zeta_n'(t-s)\zg_1\phi_n(x)\right ] f(x,s)\ZD s\,.
\end{array}
\right.
\end{equation}  

  Let the target $(\xi,\eta)\in L^2(\ZOMq)\times H^{-1}(\ZOMq)$ have the expansion 
\[
\xi(x)=\ZSUno \xi_n\phi_n(x)\,,\qquad \eta(x)=\ZSUno\left ( \eta_n\zl_n \right )\phi_n(x)\,.
\]
Then, this target is reachable  at time $T$ if and only if there exists $f\in L^2(0,T;L^2(\Gamma))$ such that   
\[
\Mo _0f= c_n = \eta_n+i\xi_n 
\] 
where $\Mo_0$ is the \emph{moment operator}
\begin{equation}\ZLA{eq:operMoment0}
\begin{array}{l}
\displaystyle \Mo_0 f=\int_\Gamma\intT E_n^{(0)}(s)\Psi_n f(x,T-s)\ZD s\ZD \Gamma \,,\\[3mm]  
 \displaystyle    \Psi_n=\frac{\zg_1\phi_n}{\zl_n}\,,\quad E_n^{(0)}(s)=\left [\zeta_n'(s)+i\zl_n \zeta_n(s)\right ] \,. 
 \end{array}
\end{equation} 
Note that the operator $\Mo_0$ takes values in $l^2$. So, we should write $\Mo_0 f=\{c_n\}$. The brace here is usually omitted. 

 The sequence $\{\Psi_n\}$ is bounded in $L^2(\partial\ZOMq)$ and it is almost normalized  if $\Gamma$ has been chosen as in Theorem~\ref{Theo:controLdue} (see~\cite[Theorem~4.4]{PandLibro}).

It is known from~\cite{LionsLibro,PandLibro} that the operator $\Mo $ is   continuous.
Our assumption is that $T$ and $\Gamma$ have been so chosen that   this operator, defined on $L^2(0,T;L^2(\Gamma))$, is surjective in $l^2(\mathbb{C})$. 
This implies (see~\cite[Sect.~3.3]{PandLibro}) that:
\begin{Lemma}\ZLA{lemmaADDE}
\begin{enumerate}
\item \ZLA{lemmaADDEitem1}
The sequence $\{e_n\}$ 
\[
e_n=\left [\zeta_n'(s)+i\zl_n \zeta_n(s)\right ]\Psi_n
\]
is a Riesz sequence in $L^2(0,T;L^2(\Gamma))$, i.e. it can be transformed to an orthonormal sequence using a linear, bounded and boundedly invertible transformation.
\item
\ZLA{lemmaADDEitem2}
The series
\[
\sum\zaa_n\Psi_n(x)\left [\zeta_n'(s)+i\zl_n\zeta_n(s)\right ]
\]
converges in $L^2(0,T;L^2(\Gamma))$ if and only if $\{\zaa_n\}\in l^2$.
\item\ZLA{lemmaADDEitem3}
If $K(t)=0$ then $
e_n(x,t)=\Psi_n(x) e^{i\zl_n t}$.
 
\end{enumerate}
\end{Lemma}
\section{\ZLA{sect:sisteVontroREG}The system with  controls of class $\Ho$}

The property $f\in \Ho$
can be written as follows:
\begin{equation}\ZLA{eq:RapprF}
\begin{array}{l}\displaystyle
f(x,t)=\intt g(x,s)\ZD s\qquad g\in \No\,,\\
\displaystyle \No=\left \{ g\in L^2(0,T;L^2(\Gamma))\,,\quad \intT g(x,s)\ZD s=0\right \}\,.
\end{array}
\end{equation}
So, 
when $f\in \Ho$  we can integrate by parts the integral in~(\ref{eqDIu})  (with $F=0 $, see~\cite{PandAMO} for the rigorous justification) and we get (using that the initial conditions are zero): 
\begin{align}
\nonumber u_f(t)&=-\ZA\intt R_-(t-s)D\ints g(r)\ZD r\,\ZD s=\\
\ZLA{eq:ondINTepartiU}&=D\intt g(r)\ZD r-\intt R_+(t-s)Dg(s)\ZD s=\hug(t)\,,\\
\ZLA{eq:ondDerivINTepartiU}u'_f(t)& =-\ZA\intt R_-(t-s)Dg(s)\ZD s=\tug(t)\,.
\end{align} 
\begin{Remark}\ZLA{Rema:sullaroBUSTrezzO}{\rm
These expressions show an interesting fact (compare also~\cite[Corollary~1.5]{EverdZUAZArch}). We see from Theorem~\ref{teo:propriESoluOnde} that  the integrals take  values respectively in $H^1_0(\ZOMq)$ and $L^2(\ZOMq)$ and
 $g\mapsto \tug(t) $ is a linear continuous map from $\No$ to $C([0,T];L^2(\ZOMq))$. Instead, $g\mapsto \hug(t)$ is a linear continuous map from $\No$ to $C([0,T]; H^{1/2}(\ZOMq))$ since ${\rm im}\,D\subseteq H^{1/2}(\ZOMq)$. We have    $\hug(T)\in H^1_0(\ZOMq)$ only 
 when $t=T$.

Of course these maps are continuous among the specified spaces but in every numerical computation we expect that the value of $T$ is affected by some error, and the fact that $u_f(t)\notin H^1_0(\ZOMq)$ if $t\neq T$ might raise some robustness issues in the numerical approximation of the steering control. This issue seemingly is still to be   studied.\zdia

}
\end{Remark}
Let us introduce

\begin{align} 
\nonumber &\No\ni g \mapsto  \Lambda_{E }^{(1)}(T)g    = \left (-\hug(T),-\tug(T)\right ) \\
\ZLA{eq:defiLambda0}=&\left (\intT R_+(T-s)Dg(s)\ZD s\,,\ \ZA\intT R_-(T-s)Dg(s)\ZD s\right ) \, .
\end{align}

Controllability with $H^1_0$-controls $f$ of the wave equation (the case $K=0$, proved in~\cite{EverdZUAZArch}) is equivalent to surjectivity of the map $ \Lambda_{E }(T)$
  from $\No$ to $H^1_0(\ZOMq)\times L^2(\ZOMq)$.

We introduce formulas~(\ref{eq:ondINTepartiU})-(\ref{eq:ondDerivINTepartiU}) in~(\ref{eq:soluVisco}) and~(\ref{eq:volteDERIVATA}). We get:

\begin{align}
\ZLA{eq:ondINTepartiW}w_f(t)&=\hug(t)+\ZA^{-1}\intt\left [ \int_0^{t-s}K(r)R_-(t-s-r)w(s)\ZD r\right] \ZD s=\hwg(t)\\
 \ZLA{eq:ondDerivINTepartiW}w'_f(t)&= \tug+ \intt \left [\int_0^{t-s}K(r) R_+(t-s-r) w(s)\ZD r\right ]\ZD s=\twg(t)\,.
\end{align}

We introduce  the operator $\Lambda_{V}^{(1)}(T)$, to be compared with the operator $\Lambda_E^{(1)}(T)$ in~(\ref{eq:defiLambda0}):
\[
\Lambda^{(1)}_{V}(T)g=\left (-\hwg(T),-\twg(T)\right )\,,\qquad g\in\No\,.
\]

Controllability in $H^1_0(\ZOMq)\times L^2(\ZOMq)$ is equivalent to surjectivity of the map $\Lambda^{(1)}_{V}(T)$ from $\No$ to $H^1_0(\ZOMq)\times L^2(\ZOMq)$.

We see from here that the functions $\hwg(t)$ and $\twg(t)$ have the same properties as stated above for $\hug(t)$ and $\tug(t)$, in particular $\hwg(t)\in H^{1/2}(\ZOMq)$ and $\hwg(T)\in H^1_0(\ZOMq)$ but there is an additional difficulty: if we want to consider $\hug(T)$ we can simply ignore the contribution of $Dg(T)$. Instead, due to the Volterra stucture of Eq.~(\ref{eq:ondINTepartiW}), the term $Dg(t)$ which comes from $\hug(t)$ cannot be simply ignore when looking at the function $\hat w_g$ for $t=T$.

In spite of this, using ${\rm im}\, D\subseteq H^{1/2}(\ZOMq)={\rm dom}\, (-A)^{1/4-\ZEP/2}$ (for every $\ZEP>0$) and solving~(\ref{eq:ondINTepartiW})-(\ref{eq:ondDerivINTepartiW})   
via Picard iteration,   it is simple to 
prove\footnote{we use 
 $\left  [\cdot \right]^\perp $ 
 to denote the subspace of the annihilators in the dual space.}:

\begin{Lemma}
\ZLA{Lemma:compa} We have:
\begin{itemize}
\item $w_f(T)\in H^{1}_0(\ZOMq)$, $w_f'(T)\in L^2(\ZOMq)$;
\item The operator $\Lambda^{(1)}_{V}(T)-\Lambda^{(1)}_{E}(T)$ is compact in $H^{1}_{0}(\ZOMq)\times L^2(\ZOMq)$ and so
${\rm Im}\, 
\Lambda^{(1)}_{V}(T)$ is closed in $H^{1}_{0}(\ZOMq)\times L^2(\ZOMq)$ and 
 $\left  [{\rm Im}\, 
\Lambda^{(1)}_{V}(T)
\right ]^\perp$ is finite dimensional.
\end{itemize}
\end{Lemma}

The goal is the proof of Theorem~\ref{Theo:controLdue}, i.e. the proof  that every   element in $\left [{\rm Im}\, \Lambda^{(1)} _{V}(T)\right ]^{\perp}$
is equal zero. We prove this result by using the properties of the moment operator.

\subsection{Controllability with $\Ho$ controls}
Our point of departure is the expansion~(\ref{eq:expaWn}) and 
the representation~(\ref{RApprewN}) of
$w_n(t)$, $w_n'(t)$.  
Let
\[
K_1(t)=\intt K(s)\ZD s
\]
 When $f$ has the representation~(\ref{eq:RapprF}) we can manipulate~(\ref{RApprewN}) as follows:
 
\begin{equation}
\ZLA{eq:RApprewNDopoINTEparti}
\begin{array}{l}
\displaystyle w_n(t)=\left (-\frac{1}{\zl_n^2}\right )\int_\Gamma \zg_1\phi_n \intt \left (-\zl_n^2\zeta_n(s)\right )\int_0^{t-s} g(r)\ZD r\,\ZD s\,\ZD\Gamma=\\[5mm]
\displaystyle =\left (-\frac{1}{\zl_n^2}\right )\int_\Gamma \zg_1\phi_n \intt\left [\zeta_n''(s)-\int_0^s K(s-\nu)\zeta_n(\nu)\ZD\nu\right ]\int_0^{t-s}g(r)\ZD r\,\ZD s\,\ZD\Gamma=\\[5mm]
\displaystyle  =\frac{1}{\zl_n^2}  \int_\Gamma \zg_1\phi_n\intt g(x,t-r)\ZD\Gamma\,\ZD r-\\
\displaystyle  -\frac{1}{\zl_n^2}\int_\Gamma\zg_1\phi_n\intt g(x,t-r)\left [
\zeta_n'(r)-\intr K_1(r-\nu)\zeta_n(\nu)\ZD\nu
\right ] \ZD r\,\ZD\Gamma
\\[5mm]
\displaystyle w_n'(t)=\int_\Gamma \zg_1\phi_n\intt g(x,t-r)\zeta_n(r)\ZD r\,\ZD \Gamma\,.
\end{array}
\end{equation}
 
Let now
\[
H^1_0(\ZOMq)\ni \xi=\ZSUno\frac{\xi_n}{\zl_n}\phi_n(x)\,,\quad
L^2(\ZOMq)\ni \eta=\ZSUno\eta_n\phi_n(x)\,,\qquad \{\xi_n\}\,,\ \{\eta_n\}\in l^2 
\]
and let \[
c_n=-\xi_n+i\eta_n
\]
Of course, $\{c_n\}$ is an arbitrary element of $l^2=
 l^2(\mathbb{C})$ and our goal is the proof that the following \emph{moment problem} is solvable for every $\{c_n\}\in l^2$:
 \begin{align*}
&\int_\Gamma\intT g(T-r)E_n^{(1)}(r)\Psi_n\ZD r\,\ZD\Gamma=c_n\,,\\
& E_n^{(1)}(r)= \left ( \zeta_n'(r)-\intr K_1(r-\nu)\zeta_n(\nu)\ZD\nu  \right )+i\zl_n \zeta_n(r) \,.
 \end{align*}
 Here $g$ \emph{is not} an arbitrary $L^2$ function, i.e. this is not a moment problem in the space $L^2\left (0,T;L^2(\Gamma)\right )$; \emph{it is a moment problem in the Hilbert space $\No$.} So, it is not really $E_n^{(1)}(r)\Psi_n$ which enters this moment problem but any projection  of $E_n^{(1)}(r)\Psi_n$ on the Hilbert space $\No$: the moment problem to be studied is
  \begin{equation}\ZLA{ProbleMOMEproiettato}
\Mo_1 g=c_n 
=\int_\Gamma\intT g(T-r){\mathcal P}_{\No}\left (E^{(1)}_n(\cdot)\Psi_n\right )\ZD r\,\ZD\Gamma 
 \end{equation}
%  where $E_n(r)=E_n(x,r)$
% is
% \[
% E_n(r) = {\mathcal P}_{\No}\left ( \Psi_n \left [\left ( \zeta_n'(r)-\intr  K_1(r-\nu)\zeta_n(\nu)\ZD\nu  \right )+i\zl_n \zeta_n(r)\right ]   \right )
% \]
  where ${\mathcal P}_{\No}$ is \emph{any fixed projection on $\No$. } 
 the operator $\Mo_1$ is the \emph{moment operator} of our control problem.
  
  So, the controllability problem    can be reformulated as follows: \emph{to prove the existence of a suitable projection ${\mathcal P}_{\No}$ such that the   moment problem~(\ref{ProbleMOMEproiettato}) is solvable for every $\{c_n\}\in l^2$.} In fact, surely there exist projections for which the moment problem is not solvable: the projection $Ph=0$ for every $h$ is an example.
 
We are going to prove that the following special projection does the job:
\begin{equation}\ZLA{eq:Laproiez}
\left ( {\mathcal P}_{\No} f\right )(t)=f(t)-\frac{1}{T}\intT f(s)\ZD s\,.
\end{equation}
\begin{Remark}
{\rm
The projection ${\mathcal P}_{\No}$ is the orthogonal projection of $L^2(0,T;L^2(\Gamma))$ onto ${\No}$. In fact, for every $f\in L^2(0,T;L^2(\Gamma))$ and every $g\in {\No}$ we have
\begin{align*}
\int_\Gamma\intT\overline{g}(x,t)\left [f-{\mathcal P}_{\No}f\right ](x,t)\ZD t\ZD\Gamma=\\
=\frac{1}{T}\int_\Gamma
\left [\intT f(x,s)\ZD s\right ]\left [\intT \overline{g}(x,t)\ZD t\right ]\ZD\Gamma=0\,.\zdiaform
\end{align*}

}
\end{Remark}
Let us note that the results reported in Section~\ref{sect:sisteVontroREG} in particular show that 
 the operator $\Mo_1$ is continuous
and the image of $\Mo_1$ is closed with finite codimension (this is Lemma~\ref{Lemma:compa}).
 
So, it is sufficient that we prove that  if $\{\bar \zaa_n\}\perp {\rm im}\, \Mo_1$ then $\{\zaa_n\}=0$ i.e. we must prove that
\[
\bigl  (\, \langle    \zaa_n,\Mo_1 g\rangle=0 \qquad \forall g\in \No\, \bigr )\ \implies \ \{\zaa_n\}=0
\]
i.e. we prove that if the following equality holds then $\{\zaa_n\}=0$:
\begin{equation}\ZLA{eq:LaSerConProi}
\ZSUno \zaa_n  \mathcal{P}_{\No}\left (\Psi_n\left [
\zeta_n'\ZCD-\int_0^{\ZCD}  K(\cdot-\nu)\zeta_n(\nu)\ZD\nu  +i\zl_n\zeta_n\ZCD
\right ]\right )=0\,.
\end{equation}
 
We introduce explicitly the projection~(\ref{eq:Laproiez}) and we see that we must prove $\{\zaa_n\}=0$ when the following equality holds:
\begin{align}
\nonumber
&\ZSUno \zaa_n\Psi_n  \left [
\zeta_n'(t)-\int_0^t  K_1(t-\nu)\zeta_n(\nu)\ZD\nu +i\zl_n\zeta_n(t)
\right ]=\\
\ZLA{eq:IlprobleFINA}
&= \frac{1}{T}
\ZSUno \zaa_n\Psi_n  \left [
\zeta_n(T)-\int_0^T\intt  K_1(t-\nu)\zeta_n(\nu)\ZD\nu\ZD t +i\zl_n\intT\zeta_n(s)\ZD s
\right ]
\end{align}
 
Note that the numerical series on the right side of~(\ref{eq:IlprobleFINA}) converges since both the series~(\ref{eq:LaSerConProi}) and the series on the left of~(\ref{eq:IlprobleFINA}) converges, 
thanks to    Lemma~\ref{lemmaADDE}. The series on the right side of~(\ref{eq:IlprobleFINA}) is constant. 
This implies that also the sum of the series on the left is constant and so its derivative is equal zero.

We shall prove that the derivative can be computed termwise. Accepting this fact, the proof that $\{\zaa_n\}=0$ is simple: the termwise derivative is
\begin{align}
\nonumber0=\ZSUno\zaa_n\Psi_n\left [\zeta_n''(t)-\intt K(t-\nu)\zeta_n(\nu)\ZD\nu+i\zl_n\zeta_n'(t)\right ]=\\
\nonumber =\ZSUno\zaa_n\Psi_n\left [-\zl_n^2\zeta_n(t)+i\zl_n\zeta_n'(t)\right ]=\\
 \ZLA{eq:laserieconSECOmembCOSTA}= i \ZSUno \zl_n\zaa_n\Psi_n\left [\zeta_n'(t)+i\zl_n
  \zeta_n(t)\right ]\,.
\end{align} 
We noted (in Lemma~\ref{lemmaADDE} item~\ref{lemmaADDEitem1})  that   $L^2(\ZOMq)\times H^{-1}(\ZOMq)$-controllability with square integrable controls  of the viscoelastic system is equivalent to the fact that
$
\left \{ \Psi_n\left [\zeta_n'(t)+i\zl_n\zeta_n(t)\right ]\right \}
$ 
 is a Riesz sequence in $ L^2(0,T;L^2(\Gamma))$ and so $\{\zl_n\zaa_n\}=0$ i.e. $\{\zaa_n\}=0$. Of course we implicitly used $\{\zl_n\zaa_n\}\in l^2$, a fact we shall prove now.

In order to complete the proof we must see that it is legitimate to distribute the derivative on    the series~(\ref{eq:IlprobleFINA}) and that this  implies in particular  $\{\zl_n\zaa_n\}\in l^2$.

The fact that $\left \{\Psi_n[\zeta_n'+i\zl_n\zeta_n]\right \}$ is a Riesz sequence in $L^2(0,T;L^2(\Gamma))$ shows that we can distribute the series on the left hand side of~(\ref{eq:IlprobleFINA}), which can be written as
\begin{align} 
\nonumber & \ZSUno \zaa_n\Psi_n  \left [
\zeta_n'(t) +i\zl_n\zeta_n(t)
\right ]\\ 
\ZLA{eq:IntRfGperH1DArichiaH2}
&=
\int_0^t   K_1(s-\nu)\left [ \ZSUno \zaa_n\Psi_n\zeta_n(\nu)\right ]\ZD\nu+{\rm const} \,.
\end{align}
So, it is sufficient that we study the differentiability of the series
\[
\ZSUno \zaa_n\Psi_n Z_n(t)\,,\qquad Z_n(t)= \left [
\zeta_n'(t) +i\zl_n\zeta_n(t)
\right ]
\,.
\]

Using the definition~(\ref{eq:diZETAn}) we see that
\begin{align*}
\zeta_n(t)=\frac{1}{\zl_n} \sin\zl_n(t)+\frac{1}{\zl_n}\left [ \intt \ints K(s-\zt)\sin\zl_n\zt\ZD\zt\right ] \zeta_n(t-s)\ZD s\,,\\
\zeta_n'(t)=\cos\zl_n(t)+\frac{1}{\zl_n} \intt \left [\ints K(s-\zt)\sin\zl_n \zt\ZD \zt\right ] \zeta_n'(s)\ZD s
\end{align*}
and so we get the following formula for $Z_n(t)$:
\begin{equation}
\ZLA{eqDiZetaROMANO}
\left\{\begin{array}{l}
\displaystyle Z_n=E_n+\frac{1}{\zl_n} K*S_n*Z_n\\[2mm]
\displaystyle \mbox{where $*$ denotes the convolution and where we define}\\[2mm]
   S_n(t)=\sin\zl_n t\,,\quad C_n(t)=\cos\zl_n t\,,\quad E_n(t)=e^{i\zl_n t}\,.
 \end{array}\right.
\end{equation}

 Gronwall inequality shows that $\{Z_n(t)\}$ is bounded on bounded intervals.

We introduce the notations
\[
F^{(*1)}(t)=F(t)\,,\quad   F^{(*k)}=F*F^{(*(k-1))}\,.
\]
With these notation the formula for $Z_n(t)$ shows also that
\begin{align}
\nonumber Z_n&=E_n+\frac{1}{\zl_n}K*S_n*E_n+\frac{1}{\zl_n^2} K^{(*2)}*S_n^{(*2)}*Z_n\,,\\
\label{asympZn}&= E_n+\sum _{r=1}^K \frac{1}{\zl_n^r}K^{(*r)}* S_n^{(*r)}*E_n+\frac{1	 }{\zl_n^{k+1}}P_{K,n}(t)\,.
\end{align}
 
The functions   $P_{n,K}(t)$ and $P'_{n,K}(t)$ are bounded on bounded intervals,
\[ 
|P_{n,K}(t)|<M_K\,,\quad |P_{n,K}(t)|<M_K
 \]
 where $ M_K $ does not depend on $ n $ and $ t\in [0,T] $.
 
  So, using the fact that $\{\Psi_n\}$ is bounded in $L^2(\Gamma)$ and~(\ref{eq:stimAUTOV}) (see~\cite[Lemma~4.4]{PandLibro})  
  we see that 
\[
\ZSUno \zaa_n\Psi_n\frac{1}{\zl_n^K}M_{n,K}(t)
\]
is of class $C^1$ (and termwise differentiable) when $K$ is large enough.
We fix an index $K$ with  this property and we consider the series of each one of the terms in~(\ref{asympZn}) for which $r\geq 1$:
\begin{equation}\ZLA{eq:LeDuESerIe}
\ZSUno\zaa_n\Psi_n\frac{1}{\zl_n}K*S_n*E_n\,,\qquad \ZSUno\zaa_n\Psi_n\frac{1}{\zl_n^r}K^{(*r)}*S_n^{(*r)}*E_n\,.
\end{equation}
$L^2$-convergence of the series is clear. We prove that they converge to $H^1$-functions.

We consider the first series. We compute the convolution $(S_n*E_n)$ and we see that this series is equal to
\begin{align*}
&\frac{1}{2i}\ZSUno\zaa_n\Psi_n \frac{1}{\zl_n}\intt K(s)\left [
(t-s)e^{i\zl_n(t-s)}-\frac{1}{\zl_n} \sin\zl_n(t-s) 
\right ]\ZD s
%=\\
%&=\frac{1}{2i} \intt K(t-s)\left [
%%%%%%%%%%%%
% \ZSUno\zaa_n\Psi_n \frac{1}{\zl_n} \left (
%se^{i\zl_ns}-\frac{1}{\zl_n} \sin\zl_ns 
%\right )\right ]\ZD s
\,.
\end{align*}
The series inside the integral converges in $ L^2(0,T) $ thanks to item~\ref{lemmaADDEitem3} in Lemma~\ref{lemmaADDE}. 
Its termwise derivative is:
\begin{align*}
&\frac{1}{2i}\ZSUno \zaa_n\Psi_n\intt K(s)\left [
\frac{1}{\zl_n} e^{i\zl_n}(t-s)+i(t-s)e^{i\zl_n (t-s)}-\cos\zl_n (t-s)
\right ]\ZD s=\\
&\frac{1}{2i}\intt K(s)\left [ 
\ZSUno \zaa_n\Psi_n \left (
\frac{1}{\zl_n} e^{i\zl_n}(t-s)+i(t-s)e^{i\zl_n (t-s)}-\cos\zl_n (t-s)\right )
\right ]\ZD s\,.
\end{align*}
This series is $L^2$-convergent thanks to Lemma~\ref{lemmaADDE}.

A similar argument holds for every $ r\leq K $. So, $ \ZSUno \zaa_n\Psi_nE_n $ converges in $ H^1(0,T) $ for every $ T $, in particular $ T>T_0 $.
  
From~\cite[Lemma~3.4]{PandLibro} and the appendix, we see that $\zaa_n=\ZDE_n/\zl_n$, $\{\ZDE_n\}\in l^2$, and that the derivative of the series can be computed termwise.

\emph{This ends the proof of item~\ref{item1delTeoINh10}  in Theorem~\ref{TeoINh10}.}

 \begin{Remark}
 {\rm
 This prrof holds in particular if $K=0$, and gives an alternative proof to the result in~\cite{EverdZUAZArch}. An important additional property in~\cite{EverdZUAZArch} is that a smooth steering control solves an optimization problem (and that its essential support is relatively compact in    $[0,T]\times\Gamma$).\zdia
 }
 \end{Remark}

\section{\ZLA{Sect:MoreREGULARcontrols}When the control is smoother}

In this section we prove   item~\ref{item2delTeoINh10} and~\ref{item3delTeoINh10} in Theorem~\ref{TeoINh10}.

We note that $f\in H^2_0(0,T;L^2(\Gamma))$ if and only if
\begin{align*}
&f(t)=\intt (t-s)g(s)\ZD s\,,\\
 & g\in \Nod=\left \{
g\in L^2(0,T;L^2(\Gamma))\,,\ \intT g(s)\ZD s=0\,,\quad \intT (T-r)g(r)\ZD r=0
\right \}\,.
\end{align*}

An analogous representation    holds if $f\in H^k_0(0,T;L^2(\Gamma))$.

Using these characterizations, we integrate by parts formulas~(\ref{eqDIu}) (with $F=0$) and~(\ref{eq:diuPRIMO}) (with zero initial conditions) and when $f\in H^2_0$ we find
\begin{equation}\ZLA{eq:LauEuprimoCONfH2}
\left\{\begin{array}{ll}
u(t)=D\intt (t-r)g(r)-\ZA^{-1} \intt R_-(t-s)Dg(s)\ZD s=\hat u_g(t)\,,\\
u'(t)=D\intt g(r)\ZD r-\intt R_+(t-s)Dg(r)\ZD r=\tilde u_g(t)\,.
\end{array}\right.
\end{equation}
This is similar to~(\ref{eq:ondINTepartiU}) and~(\ref{eq:ondDerivINTepartiU}) (and now Remark~\ref{Rema:sullaroBUSTrezzO} applies both to $u$ and $u'$).
Let
\[
%\left\{\begin{array}{ll}
L_-(t)w=\intt K(r)R_-(t-r)w\ZD r\,,\qquad
L_+(t)w=\intt K(r)R_+(t-r)w\ZD r\,.
%\end{array}
%\right.
\]
We have, from~(\ref{eq:soluVisco}) and~(\ref{eq:volteDERIVATA}),
\begin{equation}
\ZLA{eq:LauEwprimoCONfH2}
\left\{\begin{array}{ll}
w_f(t)=\hat u_g(t)+\ZA^{-1}\intt L_-(t-s)w(s)\ZD s\,,\\
w_f'(t)=\tilde u_g(t)+\intt L_+(t-s) w(s)\ZD s\,.
\end{array}\right.
\end{equation}
The proof of item~\ref{item2delTeoINh10}  in Theorem~\ref{TeoINh10} consist in 
two parts: first we prove the regularity property  $\left (w_f(T),w_f'(T)\right )\in {\rm dom}\,\ZA^2\times {\rm dom}\,\ZA$ and then we prove that $f\mapsto\left (w_f(T),w_f'(T)\right )\in  {\rm dom}\,\ZA^2\times {\rm dom}\,\ZA$ is surjective in this space.

The proof of item~\ref{item3delTeoINh10}  in Theorem~\ref{TeoINh10} is the proof that the corresponding regularity does not hold, i.e. that {due to the memory, there exist functions $f\in H^k(0,T];L^2(\Gamma))$ such that $\left (w_f(T), w'_f(T)\right )\notin {\rm dom}\,\ZA^k\times {\rm dom}\,\ZA^{k-1}$.}  
 
We proceed as follows: we first examine the regularity issue (i.e. the positive result, when $f\in H^2_0(0,T;L^2(\Gamma))$ and the lack of regularity if $f$ is smoother) in subsection~\ref{Psubs:PiuregoMENOrego}. In the subsection~\ref{H2Controll} we prove $\left ({\rm dom}\,\ZA^2\times{\rm dom}\,\ZA\right )$-controllability when $f\in H^2_0(0,T;L^2(\Gamma))$.
 
\subsection{\ZLA{Psubs:PiuregoMENOrego}Regularity and lack of regularity}
We consider $w_f(t)$ in~(\ref{eq:LauEwprimoCONfH2}). A step of Picard iteration gives
\begin{equation}\ZLA{eq:diwinH2}
w_f(t)=\hat u_g(t)+\ZA^{-1}\intt L_-(t-s)\hat u_g(s)\ZD s+A^{-1}  \left (L_-^{(*2)}*w\right) (t)\,.
\end{equation}
We know from~\cite{EverdZUAZArch} that $\hat u_g(T)\in H^2_0(\ZOMq)={\rm dom}\,A={\rm dom}\,\ZA^2$ (also seen from~(\ref{eq:LauEuprimoCONfH2})). As we noted, $\hat u_g(t)$ has this regularity for $t=T$ but not for $t\in(0,T)$. Instead, we prove that $w_f(t)-\hat u_g(t)\in {\rm dom}\,\ZA^2$ for every $t\in [0,T]$. This is clear for the third addendum on the right side of~(\ref{eq:diwinH2}).  We examine the second term, which is
\[
 \ZA^{-1}\intt L_-(t-s)D\ints(s-r) g(r)\ZD r\,\ZD s-A^{-1}  \left (L_-* R_-*g\right )(t) \,.
\]
The second addendum takes values in ${\rm dom}\, A$. We integrate by parts the first addendum and we get that it is equal to
\begin{align*}
 & \ZA^{-1}\intt L_-(t-s)D\ints(s-r) g(r)\ZD r\,\ZD s\\ 
 &=
 \ZA^{-1}\intt K(t-r)\intr R_-(r-s)D\ints (s-\nu) g(\nu)\ZD\nu\,\ZD s\,\ZD r\\
&= -A^{-1}\intt K(t-r)\left [
D\intr (r-\nu)g(\nu)\ZD\nu\right.\\
&\left.-\intr R_+(r-s)Dg(s)\ZD s
\right ]\ZD r\in {\rm dom }\,A\,.
\end{align*}

The fact that $w_f'(T)\in {\rm dom}\,\ZA$ follows from the representation of $w_f'(T)$ in terms of $w_f(t)$ in the second line of~(\ref{eq:LauEwprimoCONfH2}). This fact shows  that the integral even takes values in ${\rm dom}\,\ZA^2={\rm dom}\,A$ for every $t\in [0,T]$ while the first addendum $\tilde u_g(T)\in {\rm dom}\,\ZA$ from~\cite{EverdZUAZArch}    (also seen from~(\ref{eq:LauEuprimoCONfH2})).

 In conclusion, we proved that 
 \[
f\in H^2_0(0,T;L^2(\Gamma))\ \implies \left (w_f(T),w'_f(T)\right )\in {\rm dom}\,\ZA^2\times {\rm dom}\,\ZA\,. 
 \]
 Now we prove that this result cannot be improved when $f$ is smoother. It is sufficient that we show
 \[
f\in H^3(0,T;L^2(\Gamma))\centernot\implies w_f(T)\in {\rm dom}\,\ZA^3\,. 
 \]
 Note that $f\in H^3(0,T;L^2(\Gamma))$ when
\[
f(t)=\intt (t-s)^2g(s)\ZD s\,,\qquad \intT (T-s)^kg(s)\ZD s=0 \quad k=0\,,\ 1\,,\ 2\,.
\]
We use this representation of $f$ and we integrate by parts the integral in~(\ref{eqDIu}) with $F=0$. We get
\begin{align*}
u_f(t)=
D\intt (t-r)^2g(r)\ZD r+2A^{-1}D\intt g(r)\ZD r\\
+\ZA^{-3}\left [ -2\ZA\intt R_+(t-s)Dg(s)\ZD s\right ] \,.
\end{align*}
 The last addendum belongs to ${\rm dom}\,\ZA^3$. In the following computation we write $G(t)$ for any term which takes values in ${\rm dom}\,\ZA^3$, not the same at every occurrence. So,
 \[
 u_f(t)=\intt P(t-r)g(r)\ZD r+ G(t)\,,\qquad P(t)=t^2D+2A^{-1}D\,.
 \]
 
We replace in the expression of $w_f(t)$ in~(\ref{eq:soluVisco}). Two steps of Picard iteration
gives  the following representation for $w(t)$:
\begin{align*}
w(t)=\intt P(t-\nu)Dg(\nu)\ZD\nu+\ZA^{-1}\intt L_-(t-r)\intr(r-\nu)^2Dg(\nu)\ZD\nu\,\ZD r\\
 +\ZA^{-2}\intt  L_-(t-r)\intr L_-(r-s)\ints (s-\nu)^2Dg(\nu)\ZD\nu\, \ZD s\,\ZD r+G(t)\,.
\end{align*}
The integral in the second line can be integrated by parts  so that the second line takes values in ${\rm dom}\,\ZA^3$. The first addendum is zero for $t=T$. So, we must study the regularity of
\begin{align*}
\ZA^{-1}\intt L_-(t-r)\intr(r-\nu)^2Dg(\nu)\ZD\nu\,\ZD r\\
=-\ZA^{-2}\intt K(t-s)\ints \frac{\ZD}{\ZD r} R_+(s-r)\intr (r-\nu)^2Dg(\nu)\ZD\nu\,\ZD r\,\ZD s\\
=-\ZA^{-2} \intt K(t-s)\ints (s-\nu)^2Dg(\nu)\ZD\nu\,\ZD s\\
+2\ZA^{-2}\intt K(t-s)\ints R_+(s-r)\intr (r-\nu)Dg(\nu)\ZD\nu\, \ZD r\,\ZD s\,.
\end{align*}
The last integral can be integrated by parts again and subsumed in the term $G(t)$. Instead, 
 
 \[
\ZA^{-2}D\intT (T-\nu)^2\int_0^\nu K(s) g(\nu-s)\ZD s\,\ZD \nu \notin{\rm dom}\, \ZA^3
 \]
 since 
 \[
 D\intT (T-\nu)^2\int_0^\nu K(s) g(\nu-s)\ZD s\,\ZD \nu \notin{\rm dom}\, \ZA 
 \]
 as it is seen for example when
$g(x,\nu)=g_0(x)g_1(\nu)$ and $g_1$ such that
\[
\intT (T-\nu)^2\int_0^\nu K(s) g_1(\nu-s)\ZD s\,\ZD \nu\neq 0\,.
\]
 As a specific example, when $\ZOMq=(0,1)$ (and $\Gamma=\{0\}$)  then $Df=(1-x)f\not\in H^1_0(0,1)$ unless $f=0$.

Similar arguments hold   if $f\in H^k(0,T;L^2(\Gamma))$ and  $k>3$.

\subsection{\ZLA{H2Controll}${\rm dom}\, \ZA^2\times {\rm dom}\,\ZA$-controllability when $f\in H^2_0(0,T;L^2(\Gamma))$}
This part of the proof is similar to that in the case $f\in H^1_0(0,T;L^2(\Gamma))$ and it is only sketched.

A simple examination of formulas~(\ref{eq:LauEwprimoCONfH2}) and (\ref{eq:diwinH2}) shows  that the map $f\mapsto \left (w_f(T)-u_f(T),w'_f(T)-u'_f(T)\right )$ from $H^2_0(0,T;L^2(\Gamma))$ to ${\rm dom}\,\ZA^2\times {\rm dom}\,\ZA $ is compact. Hence we must prove
\[
%\left [R_{2,V}(T)\right ]^\perp=0\,,\qquad
  \left [\, \left \{ (w_f(T),w'_f(T))\,,\ f\in L^2(0,T;L^2(\Gamma))\right \} \, \right ]^\perp=\{0\}
\]
(the orthogonal is respect to ${\rm dom}\,\ZA^2\times {\rm dom}\,\ZA =\left (H^2(\ZOMq)\cap H^1_0(\ZOMq)\right )\times H^1_0(\ZOMq)$).

We use again the formulas for $w_n(t)$ and $w_n'(t)$ in~(\ref{eq:RApprewNDopoINTEparti}) where $g(t)=f'(t)$ has now to be replaced by $\intt g(s)\ZD s$. We see that
\begin{align*}
w_n(T)=\frac{1}{\zl_n}\int_\Gamma\intT g(T-r)\left [
\Psi_n\left (
\zeta_n(r)+\intr K_2(r-s)\zeta_n(s)\ZD s
\right )
\right ]\ZD r\,\ZD\Gamma\,,\\
w_n'(T)=\frac{1}{\zl_n}\int_\Gamma\intT g(x,T-r)\left [\Psi_n\left (-\zeta_n'(r)+\intr K_1(r-s)\zeta_n(s)\ZD s
\right )\right ]\ZD r\ZD\Gamma
\end{align*}
where
\[
K_2(t)=\intt K_1(s)\ZD s=\intt (t-s)K(s)\ZD s\,.
\]
We want to reach
\[
\xi(x)=\ZSUno\frac{\xi_n}{\zl_n^2}\phi_n(x)\,,\quad \eta(x)=\ZSUno \frac{\eta_n}{\zl_n}\phi_n(x)\,,\quad \{\xi_n\}\in l^2\,,\ \{\eta_n\}\in l^2\,.
\]
So, we must solve the moment problem in $\Nod$
\[
\int_\Gamma\intT g(T-r)P_2\left (\Psi_n E^{(2)}_n\right )\ZD r\,\ZD\Gamma =c_n\,,\qquad \{c_n\}=\{-\eta_n+i\xi_n\}\in l^2
\]
where $P_2$ is the orthogonal projection on $\Nod$:
\begin{equation}\ZLA{eq:defiPdue}
\begin{array}{l}
\left (P_2f\right )(x,t)=f(x,t)-sA_f-   B_f \,,\\
A_f=\frac{12}{T^3}\intT \left (s-\frac{T}{2}\right )f(x,s)\ZD s\,,\quad B_f=\frac{1}{T}\intT f(s)\ZD s-\frac{1}{2}T A_f
\end{array}
\end{equation}
 
and
\[
E_n^{(2)}(r)=\left (
\zeta_n'(r)-\intr K_1(r-s) \zeta_n(s)\ZD s
\right )+i
\left (
\zl_n \zeta_n(r)+\zl_n \intr K_2(r-s)\zeta_n(s)\ZD s\,.
\right )
\]
 
%QUI

Proceeding as in the case $f\in H^1_0(0,T;L^2(\Gamma))$ we see that we must prove $\{\zaa_n\}=0$ when the following equality holds:
\begin{align}
\nonumber \ZSUno \zaa_n\Psi_n\left( \zeta_n'(r)+i\zl_n\zeta_n(r)\right )\\
\nonumber =\ZSUno\zaa_n\Psi_n\left [\intr K_1(r-s)\zeta_n(s)\ZD s+i\intr K_2(r-s)\zl_n\zeta_n(s)\ZD s\right ]\\
\ZLA{EquaProieNo2MoME}+\ZSUno \zaa_n\Psi_n\left [sA_n+B_n\right ]
\end{align}
where $A_n$ and $B_n$ are as in~(\ref{eq:defiPdue}) with $f$ replaced by
\[
\left (\zeta_n'-K_1*\zeta_n\right )+i\zl_n\left ( \zeta_n+K_2*\zeta_n\right )\,.
\]
In fact, it is legitimate to distribute the series on the sum since   each one of these series converge because $\left\{\Psi_n\left (\zeta_n'(r)+i\zl_n\zeta_n(r)\right )\right\}$ is a Riesz sequence in $L^2(0,T;L^2(\Gamma))$
 (see Lemma~\ref{lemmaADDE} item~\ref{lemmaADDEitem1}).

 Equality~(\ref{EquaProieNo2MoME})  is similar to~(\ref{eq:IlprobleFINA}), with the right hand side of class $C^2$. So we get
 \[
\zaa_n=\frac{\beta_n}{\zl_n}\,,\qquad \{\beta_n\}\in l^2 
 \]
 and   we can compute   the termwise derivative of the series. Computing the derivative and noting that $K_1'=K$, $K_2'=K_1$ we get
 \[
\ZSUno \beta_n \Psi_n \left [-\zl_n \zeta_n(r)+i\zeta_n'(r)\right ]= \ZSUno \beta_n \Psi_n\intr K_1(r-s) \zeta_n(s)\ZD s-i\ZSUno\zaa_n\Psi_nA_n\,. 
 \]
 This is the same as~(\ref{eq:IntRfGperH1DArichiaH2}) and so we get $\{\beta_n\}=0$, i.e. $\{\zaa_n\}=0$ as wanted.

\section*{Appendix}

In this appendix we prove the following simple result, which however is crucial in the proof of our theorem (see also~\cite[p.~323]{GohbergKrein}).

\begin{Lemma}
Let $\mathbb{J}$ be a denumerable set and let  the sequence $\{e_n\}_{n\in \mathbb{J}}$  in a Hilbert space $H$ have the following properties:
\begin{enumerate}
\item\ZLA{itemref1} if $\{\zaa_n\}\in l^2$ then $\sum \zaa_n e_n$ converges in $H$;
\item\ZLA{itemref2} $\left \{\langle f,e_n\rangle\right \}\in l^2$ for every $f\in H$;
\item\ZLA{itemref3}  the  subspace $M=\left \{ \left \{\langle f,e_n\rangle\right \}\,,\ f\in H\right \}$  is closed  and its  codimension is finite, equal to $k$. 
\end{enumerate}
Under these conditions, there exists a set $K\subseteq \mathbb{J}$ of precisely $k$ indices
such that $\{e_n\}_{n\notin K}$ is a Riesz sequence.
\end{Lemma}
\zProof
Let $\Mo$ be the operator 
\[
\Mo f= \left \{\langle f,e_n\rangle\right \}\quad H\mapsto l^2\,.
\]
It is known that $\Mo$ is a closed operator (see~\cite{AvdoIVanovbook} and~\cite[Theorem~3.1]{PandLibro})  and Assumption~\ref{itemref2} shows that its domain is closed, equal to $H$. Hence it is a continuous operator.

$M={\rm im}\,\Mo$ and Assumption~\ref{itemref3} shows that there exist  $k$ (and not more) linearly independent sequences $\left \{ \zaa_n^i\right \}_{n\in\mathbb{J} }$, $i=1\,,\dots\, k$
such that
\begin{equation}\ZLA{eq:dellaDipeLINE}
\sum _{n\in \mathbb{J}} \zaa_n^i \langle f,e_n\rangle=0\qquad \forall f\in H\,.
\end{equation}

Note that the assumptions of this lemma does not depend on the order of the elements $e_n$. If we exchange the order of two elements $e_{n_1}$ ed $e_{n_2}$ then the corresponding elements $\zaa_{n_1}^i$ and $\zaa_{n_2}^{i}$  are exchanged for every $i=1\,,\dots,k$. So we can assume $\zaa_i^{i}\neq 0$ and without restriction $\zaa_i^{i}=1$. Hence, every $e_i$, $i=1\,,\dots\,, k$ is a linear combination of the elements $e_{n}$, $ n\in\mathbb{J}\setminus \{ 1\,,\dots\,,\ k\}$.

The operator
\[
\Mo' f= \left \{\langle f,e_n\rangle_{n\in\mathbb{J}\setminus K}\right \}\quad H\mapsto 	l^2(\mathbb{J}\setminus K) 
\]
has dense image in $l^2(\mathbb{J}\setminus K) $ otherwise we can find $\gamma$ orthogonal to its image, and adding $k$ entries equal to zero in front, we have an element orthogonal to ${\rm im}\,\Mo$, linearly independent of $\{\zaa_i\}$ and this is not possible.

Let $\Mo_0$ be the operator on $H$
\[
\Mo_0 f=\left \{\langle f,e_n\rangle \quad n\in K\right \}\,.
\]
Then, $\Mo$ is the direct sum $\Mo=\Mo_0\oplus \Mo'$ and   Lemma~\ref{Lemma:compa} shows that ${\rm im}\,\Mo$ is closed. Hence, 
the image of $\Mo'$ is   closed too, since ${\rm im}\,\Mo$ is closed. Consequently, $\Mo'$ is
surjective from $H$ to $l^2\left (\mathbb{J}\setminus \{ 1\,,\dots\,,\ k\}\right )$ and boundedly invertible. It follows that $\left\{ e_n\right \}_{n\in\mathbb{J}\setminus K} $ is a Riesz sequence.\zdia

\enddocument